\documentclass{article}
\usepackage[utf8]{inputenc}
\usepackage[T2A]{fontenc}
\usepackage{amsfonts,latexsym}
\usepackage{amsmath,amssymb,amsthm}

\usepackage[noadjust]{cite}

\newcommand{\iso}{\simeq}    % preferred isomorphism symbol

\newcommand{\Cc}{\mathbb{C}}
\newcommand{\N}{\mathbb{N}}

\newcommand{\oo}{\mathcal{O}}
\newcommand{\nn}{\mathcal{N}}
\newcommand{\vv}{\mathcal{V}}
\DeclareMathOperator{\SLtwo}{SL_2}
\DeclareMathOperator{\hd}{{hd}}
\DeclareMathOperator{\res}{{res}}
\DeclareMathOperator{\discr}{{discr}}

\newcommand{\supr[1]}{^{\,#1}}

\newtheorem{theorem}{Theorem}[section]
\newtheorem{proposition}[theorem]{Proposition}
\newtheorem{lemma}[theorem]{Lemma}

\theoremstyle{definition}

\theoremstyle{plain}

\long\def\myfootnote[#1]#2{\begingroup%
\def\thefootnote{\fnsymbol{footnote}}\footnote[#1]{#2}\endgroup}
\long\def\myfootnotetext[#1]#2{\begingroup%
\def\thefootnote{\fnsymbol{footnote}}\footnotetext[#1]{#2}\endgroup}
\def\fnf#1{{\setcounter{footnote}{#1}$^{\fnsymbol{footnote}}$}}

\def\notsobigoplus{{\mbox{$\bigoplus$}}}

\newcommand{\maysplit}{\discretionary{}{}{}}

\title{%Determination of the
$\SLtwo$-modules of small homological dimension}
\author{Andries E. Brouwer \& Mihaela Popoviciu\thanks{The second author
is partially supported by the Swiss National Science Foundation.}}
\date{Oct 17, 2010}

\begin{document}
\maketitle

\begin{abstract}
Let $V_n$ be the $\SLtwo$-module of binary forms of degree $n$ and let
$V = V_{n_1} \oplus \ldots \oplus V_{n_p}$.
We consider the algebra $R = \oo(V)^{\SLtwo}$ of polynomial functions
on $V$ invariant under the action of $\SLtwo$.
% Giving a complete description of $R$ is an old and difficult problem.
The measure of the intricacy of these algebras is the length
of their chains of syzygies, called {\em homological dimension} $\hd R$.
Popov gave in 1983 a classification of the cases in which $\hd R \le 10$
for a single binary form $(p=1)$ or $\hd R \le 3$ for a system of two
or more binary forms $(p>1)$.

We extend Popov's result and determine for $p=1$ the cases with
$\hd R \le 100$, and for $p>1$ those with $\hd R \le 15$.
In these cases we give a set of homogeneous parameters
and a set of generators for the algebra $R$.
\end{abstract}

\section{Introduction}

This paper has two goals.
First of all, following a suggestion by Popov, we extend the results
of Popov \cite{Popov} and determine all cases where the algebra of
simultaneous invariants of a number of binary forms has low
homological dimension.
Secondly, we determine the minimal degrees of a homogeneous system
of parameters (hsop) in these cases.
% In general, these parameters cannot be multi-homogeneous.
We also give a minimal system of generators,
confirming or correcting classical results.

Our base field is the field $\Cc$ of complex numbers. The group of all
complex 2$\times$2 matrices with determinant 1 is denoted $\SLtwo$. Let
$V_n$ be the set of binary forms (homogeneous polynomials in two
variables) of degree $n$.
If $V$ is a rational finite-dimensional $\SLtwo$-module, then there
exist $n_1,\ldots,n_p\in \N$ such that $V\iso V_{n_1}\oplus
\ldots \oplus V_{n_p}$ as $\SLtwo$-modules, and the algebra
$R := \oo(V)^{\SLtwo}$ of polynomial functions on $V$ invariant
under the action of $\SLtwo$ can be identified with the algebra
of joint invariants of $p$ binary forms of degrees $n_1,\ldots,n_p$.

The algebra $R$ is finitely generated (\cite{Hi1}), i.e. there exist
a finite number of invariants $j_1,\ldots, j_r$ of $V$ such that
$R=\Cc[j_1,\ldots,j_r]$.
% , and also Cohen-Macaulay (\cite{HoRo}).
Denote by $r$ the minimal number of generators
of $R$ and by $m$ the size of a system of parameters
of $R$ (set of algebraically independent elements $P_1,\ldots,P_m$ of $R$,
such that $R$ is integral over $\Cc[P_1,\ldots,P_m]$).
Then $m$ equals $\sum (n_i+1) - 3$ when this is positive,
%m = dim V//G
and the homological dimension $\hd R$ of $R$
equals $r-m$ (\cite[Corollary 1]{Popov}).

\medskip\noindent
Popov~\cite{Popov} classified the modules $V$ with the property that
$\hd R \le 3$,
\begin{table}
\centering
\begin{tabular}{|c|c|}
\hline
$V$ & $\hd R$ \\
\hline & \\[-3.5mm] \hline
$V_1$, $V_2$, $V_3$, $V_4$,
$2V_1$, $V_1\oplus V_2$, $2V_2$, $3V_1$ & 0 \\
\hline
$V_5$, $V_6$, & \\
$V_1\oplus V_3$, $V_1\oplus V_4$,
$V_2\oplus V_3$, $V_2\oplus V_4$, $2V_4$ & 1 \\
$2V_1\oplus V_2$, $V_1\oplus 2V_2$, $3V_2$, $4V_1$ & \\
\hline
$2V_3$ & 2 \\
\hline
$V_8$,~~ $5V_1$ & 3 \\
\hline
\end{tabular}
% \begin{center}
\caption{Popov's classification of $\SLtwo$-modules
with small $\hd R$}\label{tabpopov}
% \end{center}
\end{table}
and noticed that all of these were known classically.

In the past 25 years some progress was made and
sets of generators for $\oo(V_n)^{\SLtwo}$ were found in the cases
$n=7,9,10$ (\cite{DiLa,nonic,decimic}). The difficulty of this
problem is reflected by the large homological dimensions %(namely 25, 85, 98)
of the algebras of invariants in these cases.
For $R := \oo(V_n)^{\SLtwo}$ we have:
% This uses $$ for centering only.
% One can use \begin{center}...\end{center} but that gives
% additional whitespace. One can also use {\centering ... }
% with both { and } preceded by a blank line.
$$
\begin{tabular}{|c|c|c|c|c|c|c|c|c|c|c|}
\hline
$n$     & $1$ & $2$ & $3$ & $4$ & $5$ & $6$ & $7$  & $8$ & $9$ & $10$ \\
\hline
$\hd R$ & $0$ & $0$ & $0$ & $0$ & $1$ & $1$ & $25$ & $3$ & $85$ & $98$ \\
\hline
\end{tabular}.
$$

In this paper we
% look at other $\SLtwo$-modules with ``simple'' algebras of invariants
% and we use the method of Popov \cite{Popov} to
% extend his classification to:
extend Popov's classification to:

\begin{theorem}\label{1form}
Let $R := \oo(V_n)^{\SLtwo}$ and suppose that $\hd R \leq 100$.
Then $n \le 10$.
\end{theorem}

% This theorem answers Popov's question ``It will be interesting
% to know what is the minimal value of $\hd \oo(V_n)^{\SLtwo}$
% greater than 3''. (Namely, 25, for $n = 7$.)

\begin{theorem}\label{main}
Let $R := \oo(V)^{\SLtwo}$ where $V = V_{n_1} \oplus \ldots \oplus V_{n_p}$,
and suppose that $4 \le \hd R \le 15$.
Then we have one of the following:
$$
\begin{tabular}{|@{~}c@{~}|@{~}c@{~}|@{~}c@{~}|@{\,}c@{\,}|
*{9}{@{\hspace{1.1mm}}c@{\hspace{1.1mm}}|}*{2}{@{\,}c@{\,}|}}
\hline
$n_1,\ldots,n_p$ & $\hd$ & $m$ & hsop degrees & $r$ & $d_2$ & $d_3$ & $d_4$ & $d_5$ & $d_6$ & $d_7$ & $d_8$ & $d_9$ & $d_{10}$ & $d_{11}$ \\
\hline
$1,1,1,2$ & $4$ & $6$ & $2$ $(3 \times\!)$, $3$ $(3 \times\!)$ &
$10$ & $4$ & $6$ & & & & & & & & \\
\hline
$1,2,2,2$ & $5$ & $8$ & $2$ $(5 \times\!)$, $3$ $(3 \times\!)$ &
$13$ & $6$ & $4$ & $3$ & & & & & & & \\
\hline
$2,2,2,2$ & $5$ & $9$ & $2$ $(9 \times\!)$ &
$14$ & $10$ & $4$ & & & & & & & & \\
\hline
$1,1,2,2$ & $6$ & $7$ & $2$ $(4 \times\!)$, $3$ $(3 \times\!)$ &
$13$ & $4$ & $6$ & $3$ & & & & & & & \\
\hline
$1$ $(6 \times\!)$ & $6$ & $9$ & $2$ $(9 \times\!)$ & $15$ & $15$ & & & & & & & && \\
\hline
$1,1,3$ & $8$ & $5$ & $2$, $4$ $(4 \times\!)$ &
 $13$ & $1$ & & $8$ & $4$ & & & & & & \\
\hline
$1,2,3$ & $9$ & $6$ & $3,3,4,4,4,5$ &
 $15$ & $1$ & $3$ & $4$ & $4$ & $2$ & $1$ & & & & \\
\hline
$1,1,1,1,2$ & $9$ & $8$ & $2$ $(4 \times\!)$, $3$ $(3 \times\!)$, $6$ &
 $17$ & $7$ & $10$ & & & & & & & & \\
\hline
$1$ $(7 \times\!)$ & $10$ & $11$ & $2$ $(11 \times\!)$ & $21$ & $21$ & & & & &&&&& \\
\hline
$1,2,4$ & $11$ & $7$ &$2,2,3,3,4,5,6$ &
 $18$ & $2$ & $3$ & $2$ & $3$ & $4$ & $2$ & $1$ & $1$ & & \\
\hline
$2,2,3$ & $11$ & $7$ & $2,2,3,4,5,5,6$ &
 $18$ & $3$ & $2$ & $2$ & $4$ & $3$ & $4$ &&&& \\
\hline
$2,2,4$ & $11$ & $8$ &$2,2,2,3,3,3,4,4$ &
 $19$ & $4$ & $4$ & $5$ & $2$ & $4$ &&&&& \\
\hline
$1,2,2,2,2$ & $13$ & $11$ & $2$ $(7 \times\!)$, $3$ $(4 \times\!)$ &
 $24$ & $10$ & $8$ & $6$ & & & & & & & \\
\hline
$2,2,2,2,2$ & $13$ & $12$ & $2$ $(12 \times\!)$ &
 $25$ & $15$ & $10$ & & & & & & & & \\
\hline
$4,4,4$ & $13$ & $12$ & $2$ $(6 \times\!)$, $3$ $(6 \times\!)$ &
 $25$ & $6$ & $10$ & $6$ & $3$ & & & & & & \\
\hline
$1,1,4$ & $14$ & $6$ & $2,3,5,5,6,6$ &
 $20$ & $2$ & $1$ & & $5$ & $5$ & & & $7$ & & \\
\hline
$3,4$ & $14$ & $6$ & $2,3,4,5,6,7$ &
 $20$ & $1$ & $1$ & $1$ & $2$ & $2$ & $3$ & $3$ & $4$ & $2$ & $1$ \\
\hline
$1$ $(8 \times\!)$ & $15$ & $13$ & $2$ $(13 \times\!)$ &
 $28$ & $28$ &&&&&&&&& \\
\hline
$1,1,1,2,2$ & $15$ & $9$ & $2$ $(5 \times\!)$, $3$ $(4 \times\!)$ &
 $24$ & $6$ & $12$ & $6$ &&&&&&& \\
\hline
\end{tabular}
$$
Here $V$ has a minimal set of generators of size $r$,
with $d_i$ generators of degree $i$ $(2 \le i \le 11)$.
The size of any homogeneous system of parameters (hsop) is $m$,
and the degrees for one particular such system are as given.
The column $\hd$ gives $\hd R$.
\end{theorem}

The paper is organised as follows: In \S \ref{classical} we describe
the classical results and correct them where needed.
In \S \ref{numgens} we find a lower bound for $r$ given
the Poincar\'e series. In \S \ref{findV} we determine $V$.
In \S \ref{findgens} we describe how to find a set of generators.
A prerequisite is a homogeneous system of parameters, found in
\S \ref{hsops}. The actual generators are constructed
in \S \ref{gens}.

\section{The classical results}\label{classical}
The table below gives the classical (that is, 19th century)
results\myfootnote[2]{
In more complicated cases the classical techniques were not powerful enough
to determine the precise values of $r$---the German school found
upper bounds only, the English school claimed to find true values,
or at least lower bounds, but the former was mistaken (cf.~\cite{Hammond}),
the latter unproved.}
% and the values obtained were uncertain.}
(\cite{Bessel,Elliot,vGall26,vGall8,vGall7,vGall333,GordanS,Gordan,Gund,GrYo,
Perrin,Sinigallia,Sylv,Sylv8,Winter,Young}), possibly slightly amended.
Here $aV_s$ stands for the direct sum $\oplus_{i=1}^a V_s$ of $a$ copies of
$V_s$. See also \S\ref{mult1}.
\begin{table}[h]
\centering
\begin{tabular}{|c|c||c|c||c|c||c|c|}
% \begin{tabular}{|c|@{~}c@{~}||c|@{~}c@{~}||c|@{~}c@{~}||c|@{~}c@{~}|}
\hline
module & $r$ & module & $r$ & module & $r$ & module & $r$ \\
\hline
$V_1$ & 0 &
$2V_1$ & 1 &
$3V_1$ & 3 &
$4V_1$ & 6 \\
$V_2$ & 1 &
$V_1 \oplus V_2$ & 2 &
$3V_2$ & 7 &                        % \cite{Bessel}
$V_1 \oplus 3V_2$ & 13 \\           % \cite{GordanS}
$V_3$  & 1 &
$V_1 \oplus V_3$  & 4 &
$3V_3$ & 28 &                       % \cite{vGall333,Sinigallia}
$V_1 \oplus 3V_3$ & 97\fnf{6} \\    % vGall333: 98, Sinigallia: 97
$V_4$  & 2 &
$V_1 \oplus V_4$  & 5 &
$3V_4$ & 25 &
$V_1 \oplus 3V_4$ & 103\fnf{5} \\
$V_5$  & 4 &                        % \cite{Gordan}, p.240
$V_1 \oplus V_5$ & 23 &             % Gordan 1868, \cite{Gordan}, p.240
$4V_4$ & 80 &                       % \cite{Young}
$V_1 \oplus 4V_4$ & 305\fnf{5} \\   % \cite{Young}
$V_6$  & 5 &                        % \cite{Gordan}, p.283
$V_1 \oplus V_6$  & 26 &            % \cite{Gordan}, p.283
$V_2 \oplus V_3$ & 5 &              % \cite{Bessel}
$V_1 \oplus V_2 \oplus V_3$ & 15 \\ % Salmon and Clebsch, \cite{Gordan}, p.324
$V_7$ & 30\fnf{3} &
$V_1 \oplus V_7$ & 147\fnf{9} &
$V_2 \oplus V_4$ & 6 &              % \cite{Bessel}
$V_1 \oplus V_2 \oplus V_4$ & 18 \\ % \cite{GordanS}, Clebsch, p.214
$V_8$ & 9 & $V_1 \oplus V_8$ & 69\fnf{4} & % \cite{vGall8}
$V_2 \oplus V_5$ & 29 &             % checked
$V_1 \oplus V_2 \oplus V_5$ & 92\fnf{8} \\ % checked, cf. \cite{Winter}
$2V_2$ & 3 &                        % \cite{Bessel}
$V_1 \oplus 2V_2$ & 6 &             % \cite{GordanS}
$V_2 \oplus V_6$ & 27 &             % \cite{vGall26}
$V_1 \oplus V_2 \oplus V_6$ & 99 \\ % \cite{vGall26}
$2V_3$ & 7 &                        % \cite{Gordan}, p.333
$V_1 \oplus 2V_3$ & 26 &            % Salmon and Clebsch, \cite{Gordan}, p.333
                                    % Clebsch, p.224-225 has 28 forms.
$V_3 \oplus V_4$ & 20 &             % checked by MP
$V_1 \oplus V_3 \oplus V_4$ & 63$^*$ \\ % Gundelfinger, Sylvester, aeb
$2V_4$ & 8 &
$V_1 \oplus 2V_4$ & 28 &            % \cite{GordanS}, \cite{Sylv}, Bertini,
% $2V_1 \oplus 2V_2$ & 13 &
% $3V_1 \oplus 2V_2$ & 24 \\          % \cite{Perrin}
&&& \\
% $2V_1 \oplus V_2$ & 5 &
% $2V_1 \oplus V_3$ & 13 &            % \cite{Elliot}, p.334, \cite{Sylv}
% $2V_1 \oplus V_4$ & 20 && \\
\hline
\end{tabular}
\caption{The classical results}\label{tabclassic}
\end{table}

% footnotes do not work in tables, so have separate labels and texts
\myfootnotetext[1]{Gundelfinger found 64, Sylvester 61, it is 63.}
\myfootnotetext[4]{von Gall found 96, then 67, then 70; Sylvester 69.
See also Shioda \cite{Shi} and Bedratyuk \cite{Be08}.}
\myfootnotetext[3]{von Gall found 33, Sylvester 26, Hammond two more,
Dixmier \& Lazard 30.}
\myfootnotetext[9]{von Gall found 153, Sylvester 124,
Cr\"{o}ni \!\cite{Croeni} and Bedratyuk \!\cite{Be07} find 147.}
\myfootnotetext[6]{von Gall found 98, Sinigallia 97.
(Peano \cite{Peano} has partial results on $pV_3$, $V_1 \oplus pV_3$.)}
\myfootnotetext[8]{Winter found 94, it is 92.}
\myfootnotetext[5]{Young \cite{Young} treats $pV_4$ and $V_1 \oplus pV_4$
for all $p$.}
More generally, Gordan \cite{GordanS,Gordan} gives
for $V = pV_1$ the value $r = \binom{p}{2}$,
for $V = pV_2$ the value $r = \binom{p+1}{2} + \binom{p}{3}$,
and for $V = V_1 \oplus qV_2$ the value $r = q(q+1) + \binom{q}{3}$,
cf.~\S\ref{onetwo} below.
%
% In fact, for $V = mV_1 + nV_2$ one has $r = \binom{n}{3} +
% \binom{m+1}{2}\binom{n+1}{2} + \binom{m}{2} + \binom{n+1}{2}$,
% see~\S\ref{onetwo}. % below.
%
From the generators in case $V \oplus V_1$ one can derive those
for $V \oplus pV_1$ for all $p > 1$, cf.\,\cite[\S 55]{Clebsch},
\cite[\S 138A]{GrYo} and \S\ref{mult1} below.

\section{The number of generators}\label{numgens}
The Poincar\'e series of a graded $k$-algebra $R = \oplus_i R_i$
is defined as $P(t)=\sum_i a_i t^i$, where $a_i = \dim_k (R_i)$.
Here we consider $k = \Cc$ and $R = \oo(V)^{\SLtwo}$, where $V$ is
an $\SLtwo$-module.
% This algebra is graded by degree, so that $R = \oplus_k R_k$, where
% $R_k$ is the vector space of degree $k$ elements of $R$ (together with 0).
Formulas for the coefficients $a_i$ were already given by
Cayley and Sylvester. A closed expression for $P(t)$ as a rational
function in $t$ was given by Springer~\cite{Springer} for the case
of $V = V_n$, and by Brion~\cite{Brion} in general.
The webpages \cite{aebp} list some results of computations
due to Bedratyuk and Brouwer that we~use.

\subsection{Tamisage}
Suppose $R$ has Poincar\'e series $P(t) = \sum a_i t^i$.
(Then $a_0 = 1$ and $a_1 = 0$.)
Determine numbers $m_i$ as follows: As long as there is an $i > 0$
for which $a_i \ne 0$, find the smallest such $i$. If $a_i < 0$, stop.
Otherwise put $m_i := a_i$ and replace $P(t)$ by $P(t)(1-t^i)^{m_i}$
and repeat. Let undefined $m_i$ be zero. This is the process that
Sylvester called `tamisage'.

\medskip\noindent{\bf Sylvester's claim}.\!\!\myfootnote[1]{
``If the {\em fundamental postulate} were called into question,
this (it may be proved) would not affect the fact of the existence
of the groundforms obtained by its aid, but only the possibility
of the existence of other groundforms over and above those so obtained.
Thus my tables of groundforms could only err (were that possible, which
I do not believe it to be) in defect; and as those found by the German
method can only err in excess, it follows that, whenever the tables coincide,
both must be correct.'' (J. J. Sylvester \cite[p.\,249]{Sylv})} % 38, 309
{\it The number of generators of $R$ is at least $\sum_i m_i$.
More precisely: the number of generators of $R$ of degree $i$
is at least $m_i$.}
\medskip

So far this claim is unproved.
% Sylvester's stronger claim that his `fundamental postulate' holds
% (which implies that his tamisage gives true values rather than just
% lower bounds) is false, as was shown by Hammond \cite{Hammond}.
We use a slightly weaker bound in the below, one that has the advantage
of having an easy proof. Maintain two numbers $m_i$ and $M_i$
as lower and upper bounds for the number of generators of degree $i$
in a minimal system of generators. Also maintain upper bounds
$M_{ij}$ for the dimension of the space of degree $i$
invariants spanned by those having a factor of degree $j$ but no
factor of smaller degree, for $j \le i$.
Put $m_i = a_i - \sum_{j<i} M_{i,j}$ and
$$
M_{ij} = \min \{ a_{i-j}a_j,~
\sum_{t \ge 1,\,tj \le i} \binom{m+t-1}{t} S_{i-tj,j} \}
$$
where $m = M_j$ and $S_{0,j} = 1$ and $S_{a,j} = \sum_{k>j} M_{a,k}$
for $a > 0$. Finally put
\begin{align*}
d_1 &= \max_{j<i,\,a_j \ne 0} a_{i-j}, \\
d_2 &= \max_{j<i,\,a_j \ge 2} (2a_{i-j}-a_{i-2j}), \\
d_3 &= \max_{j<k<i,\,a_ja_k \ne 0} (a_{i-j}+a_{i-k}-a_{i-j-k}), \\
M_i &= a_i - \max \{ 0, d_1, d_2, d_3 \} .
\end{align*}
where $a_h = 0$ for $h < 0$. This satisfies all requirements.
Indeed, for $M_i$ we need to subtract from $a_i$ a lower bound
for the number of linearly independent invariants of degree $i$
that have a factor of some smaller degree. If $u \in R_j$, then
$x \mapsto ux$ is an injection of $R_{i-j}$ into $R_i$, so that
$d_1$ is such a lower bound. Now consider distinct basic invariants
$u \in R_j$ and $v \in R_k$. The images of $x \mapsto ux$ and
$y \mapsto vy$ (for $y \in R_{i-k}$)
have an intersection consisting of invariants with factor $uv$,
so that the dimension of the intersection is $a_{i-j-k}$.
This shows that also $d_2$ and $d_3$ are lower bounds.
The value given for $m_i$ is clear. Concerning $M_{ij}$,
if an invariant of degree $i$ has precisely $t$ factors
that are basic invariants of degree $j<i$, then the quotient
of degree $i-tj$ can be chosen in (at most) $S_{i-tj,j}$ ways
and the product of $t$ factors can be chosen in $\binom{m+t-1}{t}$ ways.

Now the final lower bound for the number of generators is
$
r \ge \sum_i m_i .
$

\medskip\noindent
{\bf Example.}
The Poincar\'e series $P(t)$ of $\oo(V_{12})^{\SLtwo}$ starts
\[ 1 + t^2 + t^3 +  3t^4 +  3t^5 +  8t^6 +  10t^7 +  20t^8 +
28t^9 +  52t^{10} +  73t^{11} +  127t^{12} + \ldots \]
We find
$$
\begin{tabular}{@{~}cccccccccccccccccc@{~}}
$i$ & 2 & 3 & 4 & 5 & 6 & 7 & 8 & 9 & 10 & 11 & 12 & 13 & 14 & 15 & 16 & 17 \\
\hline
$m_i$ & 1 & 1 & 2 & 2 & 4 & 5 &  7 &  9 & 14 & 12 &  9 & 0 & 0 & 0 & 0 & 0 \\
$g_i$ & 1 & 1 & 2 & 2 & 4 & 5 &  7 &  9 & 14 & 15 & 19 & 18 & 12 & 2 & 1 & 1 \\
$M_i$ & 1 & 1 & 2 & 2 & 4 & 5 & 10 & 13 & 25 & 33 & 57 & 76 &&& \\
\end{tabular}
$$
so that $r \ge 66$, $\hd R \ge 56$. The row $g_i$ gives the actual
number of generators of degree $i$ (known in this case, cf.~\cite{aebi}),
so that $r \ge 113$, $\hd R \ge 103$.

\subsection{Bounds}
In the table below, we list modules, the Poincar\'e series, and lower bounds
for $r$ and $\hd R$. In many cases, better bounds are obtained by taking
more terms.

\medskip\noindent
\begin{tabular}{|@{~}c@{~}|l@{~}|c|@{~}c@{~}|}
\hline
module & Poincar\'e series & $r \ge$ & $\hd R \ge$ \\
\hline
$V_{11}$ & $1 + 2t^4 + 13t^8 + 13t^{10} + 73t^{12} + 110t^{14} + \ldots$ &
  158 & 149 \\
% $V_{12}$ & $1 + t^2 + t^3 + 3t^4 + 3t^5 + 8t^6 + 10t^7 + 20t^8 + $&&\\
%   & ~~$28t^9 + 52t^{10} + 73t^{11} + 127t^{12} + 181t^{13} + \ldots$ &
%   66 & 56 \\
$V_{13}$ & $1 + 2t^4 + 22t^8 + 33t^{10} + 181t^{12} + 375t^{14} + \ldots$ &
  502 & 491 \\
$V_{14}$ & $1 + t^2 + 3t^4 + 10t^6 + 4t^7 + 31t^8 + 27t^9 + $ &&\\
  & ~~$97t^{10} + 110t^{11} + \ldots$ & 182 & 170 \\
$V_{15}$ & $1 + 3t^4 + t^6 + 36t^8 + 80t^{10} + 418t^{12} + \ldots$ &
  425 & 412 \\
$V_{16}$ & $1 + t^2 + t^3 + 3t^4 + 4t^5 + 13t^6 + 18t^7 + 47t^8 +$&&\\
  & ~~$84t^9 + 177t^{10} + \ldots$ & 198 & 184 \\
$V_{18}$ & $1 + t^2 + 4t^4 + t^5 + 16t^6 + 13t^7 + 71t^8 + 99t^9 +$ & 161 & 145 \\
$V_{20}$ & $1 + t^2 + t^3 + 4t^4 + 5t^5 + 20t^6 + 35t^7 + 102t^8 +$ & 123 & 105 \\
$V_{22}$ & $1 + t^2 + 4t^4 + t^5 + 24t^6 + 26t^7 + 144t^8 + \ldots$ & 164 & 144 \\
$V_{24}$ & $1 + t^2 + t^3 + 5t^4 + 7t^5 + 29t^6 + 62t^7 + 201t^8 +$ & 242 & 220 \\
$V_{28}$ & $1 + t^2 + t^3 + 5t^4 + 8t^5 + 40t^6 + 97t^7 + 365t^8 +$ & 440 & 414 \\
$V_{32}$ & $1 + t^2 + t^3 + 6t^4 + 10t^5 + 54t^6 + 153t^7 + \ldots$ & 201 & 171 \\
\hline
% \end{tabular}
% 
% \medskip\noindent
% \begin{tabular}{|c@{~}|l|c|c@{~}|}
% \hline
% module & Poincar\'e series & $r \ge$ & $\hd R \ge$ \\
% \hline
$V_2 \oplus V_8$ & $1 + 2t^2 + t^3 + 5t^4 + 5t^5 + 15t^6 + 17t^7 +$ &&\\
  & ~~$41t^8 + 54t^9 + 108t^{10} + \ldots$ & 35 & 26 \\
$V_3 \oplus V_8$ & $1 + t^2 + t^3 + 3t^4 + 4t^5 + 9t^6 + 16t^7 + 30t^8 \!+
  \ldots\!$ & 37 & 27 \\
$V_4 \oplus V_8$ & $1 + 2t^2 + 4t^3 + 8t^4 + 16t^5 + 35t^6 + 60t^7 +
  \ldots$ & 42 & 31 \\
$V_5 \oplus V_8$ & $1 + t^2 + t^3 + 3t^4 + 6t^5 + 15t^6 + 31t^7 + \ldots$ &
  43 & 31 \\
$V_6 \oplus V_8$ & $1 + 2t^2 + 2t^3 + 10t^4 + 14t^5 + 46t^6 + 82t^7 + \ldots$ &
  88 & 75 \\
% \hline
% \end{tabular}
% 
% \medskip\noindent
% \begin{tabular}{|l|l|c|c@{~}|}
% \hline
% module & Poincar\'e series & $r \ge$ & $\hd R \ge$ \\
% \hline
\hline
\end{tabular}

\medskip\noindent
\begin{table}
\begin{tabular}{|l@{~}|l@{~}|@{~}c@{~}|@{~}c@{~}|}
\hline
module & Poincar\'e series & $r \ge$ & $\hd R \ge$ \\
\hline
$V_1\oplus 2V_3$ & $1 + t^2 + 13t^4 + 26t^6 + \ldots$ & 26 & 19\\
$V_2 \oplus 2V_3$ & $1+2t^2+3t^3+9t^4+12t^5+26t^6+44t^7+\ldots $ & 26 & 18 \\
$V_1 \oplus 2V_2 \oplus V_3$ & $1+3t^2+6t^3+15t^4+30t^5+65t^6+\ldots$ & 34 & 25 \\
% $3V_2 \oplus V_3$ & $1+6t^2+4t^3+25t^4+34t^5+\ldots$ & 24 & 14 \\
\hline
$V_2 \oplus V_3 \oplus V_4$ & $1+2t^2+3t^3+7t^4+14t^5+29t^6+52t^7 +\ldots$ & 43 & 34 \\
$V_1 \oplus 2V_2 \oplus V_4$ & $1+4t^2+6t^3+18t^4+33t^5 +\ldots$ & 27 & 17 \\
$3V_2 \oplus V_4$ & $1+7t^2+8t^3+42t^4+64t^5 +\ldots$ & 37 & 26 \\
$2V_3 \oplus V_4$ & $1+2t^2+2t^3+9t^4+16t^5+37t^6+71t^7+\ldots$ & 69 & 59 \\
$V_3 \oplus 2V_4$ & $1+3t^2+4t^3+10t^4+22t^5+49t^6+96t^7+\ldots$ & 45 & 34 \\
\hline
% \end{tabular}
% 
% \medskip\noindent
% \begin{table}
% \begin{tabular}{|l|l|c|c@{~}|}
% \hline
% module & Poincar\'e series & $r \ge$ & $\hd R \ge$ \\
% \hline
$V_3\oplus V_5$ & $1+6t^4+7t^6+36 t^8 +\ldots $ & 28  & 21  \\
$V_4\oplus V_5$ & $1 + t^2 + t^3 + 2t^4 + 4t^5 + 8t^6 + 12t^7 + 22t^8 +$ &&\\
  & ~~$37t^9 + 56t^{10} + \ldots$ & 59 & 51 \\
$2V_5$ & $1 + t^2 + 7t^4 + 14t^6 + 72t^8 + 168t^{10} + \ldots$ & 105 & 96 \\
\hline
% \end{tabular}
% 
% \medskip\noindent
% \begin{tabular}{|l|l|c|c@{~}|}
% \hline
% module & Poincar\'e series & $r \ge$ & $\hd R \ge$ \\
% \hline
$V_3\oplus V_6$ & $1+t^2+t^3+3t^4+4t^5+8t^6+12 t^7+21 t^8 +\ldots$ & 24 & 16 \\
$V_4\oplus V_6$ & $1 + 2t^2 + 2t^3 + 7t^4 + 8t^5 + 24t^6 + 31t^7 +
  68t^8 +$ & 33 & 24 \\
$V_5\oplus V_6$ & $1+t^2+t^3+3t^4+5t^5+12t^6+22 t^7+\ldots$ & 31 & 21 \\
$2V_6$ & $1 + 3t^2 + 12t^4 + 6t^5 + 44t^6 + 40t^7 + 150t^8+\ldots$ & 29 & 18 \\
\hline
\end{tabular}
\caption{Bounds from the Poincar\'e series}\label{tabpoinc}
\end{table}
\vspace{-1cm}

\section{Determining $V$}\label{findV}

Consider $V=V_{n_1}\oplus \ldots \oplus V_{n_p}$ with $n_i\geq 1$ for
all $i$. Let $R := \oo(V)^{\SLtwo}$ be the algebra of invariants of $V$.
We want to determine $V$ if either $p=1$ and $\hd R \leq 100$,
or $p > 1$ and $\hd R \le 15$.

\medskip
First consider the case $p = 1$, $V = V_n$.
By \cite[Proposition~6]{Popov}, if $n$ is even and $\hd R \le 100$,
% n=2b, b odd: b=13 hd >= 121, b=11 hd >= 81
% n=2b, b even: b=18 hd >= 120, b=16 hd >= 91
then $n \le 24$ or $n \in \{28,32\}$.
By \cite[p.\,106]{Kac}, if $n$ is odd, then $r \ge p(n-2) + \phi(n-2) - 1$,
where $p()$ is the partition function and
$\phi()$ is Euler's totient function.
%
%  n r\ge hd\ge
%  3   1   0
%  5   4   1
%  7  10   5
%  9  20  13
% 11  35  26
% 13  65  54
% 15 112  99
% 17 183 168
% 19 312 295
It follows that $\hd R \ge 168$ for odd $n \ge 17$.
We know $\hd R$ for $n \le 10$ (a table was given above),
and $\hd R \ge 103$ for $n = 12$ (see the example above),
and for the remaining values we found $\hd R \ge 105$
in Table \ref{tabpoinc}. This proves Theorem \ref{1form}.

\medskip
Now consider the case $p > 1$ and assume $\hd R \le 15$.
By the monotony theorem ~\cite[Theorem 2b]{Popov} we have
if $V = W \oplus W'$, then
$\hd R \ge \hd \oo(W)^{\SLtwo} + \hd \oo(W')^{\SLtwo}$.
Therefore, all $n_i$ belong to $\{1,2,3,4,5,6,8\}$,
and direct summands $W$ have $\hd \oo(W)^{\SLtwo} \le 15$.

\medskip
If all $n_i$ are either 1 or 2, so that $V = mV_1 \oplus nV_2$,
then we have the explicit formula $r = \binom{n}{3} +
\binom{m+1}{2}\binom{n+1}{2} + \binom{m}{2} + \binom{n+1}{2}$
(see \S\ref{onetwo} below), and $\hd R = r - (3n+2m-3)$ for $m+n > 1$.
A table of $\hd R$ shows that $V$ is as claimed.
$$
\begin{tabular}{c|cccccccccc}
$n \backslash m$ & 0 & 1 & 2 & 3 & 4 & 5 & 6 & 7 & 8 & 9 \\
\hline
0 & 0 & 0 & 0 & 0 & 1 & 3 & 6 & 10 & 15 & 21 \\
1 & 0 & 0 & 1 & 4 & 9 & 16 & 25 && \\
2 & 0 & 1 & 6 & 15 & 28 &&&& \\
3 & 1 & 5 & 16 & 34 &&&&& \\
4 & 5 & 13 & 32 &&&&&& \\
5 & 13 & 26 &&&&&& \\
6 & 26 &&&&&&& \\
\end{tabular}
$$

Now investigate the remaining possibilities.

\medskip\noindent
{\em Case 1: One of the $n_i$, say $n_1$, is equal to 8.}

\noindent
By the discussion in the proof of \cite[Theorem 4]{Popov},
$\hd R \geq \frac{(b-1)(b-2)}{2}$,
where $b=3+[\frac{n_2+1}{2}]+\ldots +[\frac{n_p+1}{2}]$.
If $b\geq 8$, then $\hd R \geq 21$.
We check the cases with $b \le 7$. By monotony it suffices
to look at $V_m \oplus V_8$ for $m=1,2,3,4,5,6$.
If $V = V_1 \oplus V_8$ then $R$ is the algebra of covariants of $V_8$,
generated by 69 elements (\cite{Be08}), and $\hd R = 61$.
In the other cases $\hd R \ge 26$ by Table \ref{tabpoinc}.

\medskip\noindent
{\em Case 2: One of the $n_i$, say $n_1$, is equal to 2.}

\noindent
By the discussion in the proof of \cite[Theorem 4]{Popov},
$ \hd R \geq (c-1)^2$, where $c=[\frac{n_2+1}{2}]+\ldots +[\frac{n_p+1}{2}]$.
Since $n_i > 2$ for some $i$, we have $c\geq 2$.
We have $\hd R \ge 16$ for $c \geq 5$.
We therefore check the cases $c\in \{2,3,4\}$.

If $c=2$, then $V$ is $V_2 \oplus V_3$ or $V_2 \oplus V_4$ and $\hd R = 1$.

If $c=3$, then $V$ is one of
$V_1 \oplus V_2 \oplus V_3$, $V_1 \oplus V_2 \oplus V_4$,
$2V_2 \oplus V_3$, $2V_2 \oplus V_4$, $V_2 \oplus V_5$ or $V_2 \oplus V_6$.
In these six cases one has $\hd R = 9, 11, 11, 11, 23, 20$, respectively.
% For the first two cases, cf.~\cite{GrYo}, p.\,165,\,168.
(For the first two and last two cases, see Table \ref{tabclassic}.
For the other two, see \S\S\ref{V223},\,\ref{V224}.)

If $c=4$, then by monotony and the above $V$ does not have a direct
summand $V_5$ or $V_6$, so that $V$ is one of
$V_2\oplus 2V_3$, $V_2\oplus V_3 \oplus V_4$, $V_2\oplus 2V_4$,
$2V_1\oplus V_2 \oplus V_3$, $2V_1\oplus V_2 \oplus V_4$,
$V_1\oplus 2V_2 \oplus V_3$, $V_1\oplus 2V_2 \oplus V_4$,
$3V_2 \oplus V_3$, $3V_2 \oplus V_4$.
If $V$ is $2V_1\oplus V_2 \oplus V_3$ or $2V_1\oplus V_2 \oplus V_4$,
then $\hd R = 27$ or $48$ by Table \ref{clmult1}.
Explicit generation of invariants for $V_2\oplus 2V_4$ and $3V_2 \oplus V_3$
shows that $r \ge 29$, $49$ so that $\hd R \ge 19$, $39$ in these cases.
By Table \ref{tabpoinc} $\hd R \ge 17$ in the remaining five cases.

\medskip\noindent
{\em Case 3: All of the $n_i$ equal 1, 3, 4, 5 or 6.}

\noindent
If $V$ is $V_1 \oplus V_3$, $V_1 \oplus V_4$, $2V_3$, $V_3\oplus V_4$,
$2V_4$, $V_1\oplus V_5$, or $V_1 \oplus V_6$, then $\hd R$ equals
1, 1, 2, 14, 1, 18, 20, respectively, by Table \ref{tabclassic}.
If $V$ is $V_3\oplus V_5$, $V_4\oplus V_5$, $2 V_5$,
$V_3\oplus V_6$, $V_4\oplus V_6$, $V_5\oplus V_6$, $2V_6$,
$2V_3 \oplus V_4$ or $V_3 \oplus 2V_4$,
then $\hd R \ge 16$ by Table \ref{tabpoinc}.
If $V$ is $2V_1\oplus V_3$, $2V_1\oplus V_4$, $V_1\oplus 2V_3$,
$V_1 \oplus V_3 \oplus V_4$, $V_1\oplus 2V_4$,  $3V_3$,
$3V_4$, $3V_1 \oplus V_3$, $3V_1 \oplus V_4$, $4V_4$, 
then $\hd R$ equals 8, 14, 19, 55, 19, 19, 13, 23, 55, 63,
respectively, by Tables \ref{tabclassic} and \ref{clmult1}.
By monotony we are done.

\medskip\noindent
This finishes the determination of the $V$ with $\hd R \le 15$.

\section{Finding generators}\label{findgens}

Let $V$ be an $\SLtwo$-module, and $R = \oo(V)^{\SLtwo}$ its algebra
of invariants. Finding a minimal set of generators of $R$ is routine,
only requiring computational power, if a good upper bound for the
maximum degree of these generators is known.
Details for $V_9$ and $V_{10}$ were given in \cite{nonic,decimic}.
(The cases considered here are much smaller.)
An upper bound for the maximum degree of a basic generator follows
from the Poincar\'e series when a hsop (homogeneous system of parameters),
or at least the set of degrees of a hsop, is known.

\subsection{Hilbert's Criterion}
One way of computing a system of parameters of $R$ is
finding equations for the {\em nullcone} of $V$. The nullcone of $V$,
denoted $\nn (V)$, is the set of all elements of $V$ on which
all invariants vanish.
One shows that $\nn (V_{n_1}\oplus \ldots \oplus V_{n_p})$ is
the set of all $(f_1,\ldots,f_p)\in V_{n_1}\oplus \ldots \oplus V_{n_p}$
such that $f_1,\ldots,f_p$ have a common root of multiplicity
$>\frac{1}{2}n_i$ in $f_i$ for all $i=1,\ldots,p$. This is a consequence of
the Hilbert-Mumford criterion. % (see \cite[III.2.4.]{Kr} for details). 
Let $\vv(J)$ stand for the vanishing locus of $J$.

\begin{proposition} \label{hilb} {\rm (Hilbert \cite{Hi2})}
Let $V = V_{n_1}\oplus \ldots \oplus V_{n_p}$, and
$R = \oo(V)^{\SLtwo}$, and $m = n_1+\ldots+n_p+p-3 > 0$.
A set $P_1,\ldots, P_m$ of homogeneous elements of $R$
is a system of parameters of $R$ if and only if
$\vv(P_1,\ldots, P_m)=\nn(V)$.
\end{proposition}

\subsection{Dixmier's Criterion}
Since we do not actually need the hsop but only the degrees,
the following is often easier to apply than Hilbert's Criterion.

\begin{proposition} {\rm (Dixmier \cite{Di1})}
Let $G$ be a reductive group over $\Cc$, with a rational representation
in a vector space $V$ of finite dimension over $\Cc$. Let $\oo(V)$ be the
algebra of complex polynomials on $V$, $R := \oo(V)^G$ the subalgebra of
$G$-invariants, and $R_d$ the subset of homogeneous polynomials
of degree $d$ in $R$. Let $X$ be the affine variety such that
$\Cc[X] = R$. Let $m = \dim X$. Let $(d_1,\ldots,d_m)$
be a sequence of positive integers. Assume that for each subsequence
$(j_1,\ldots,j_p)$ of $(d_1,\ldots,d_m)$ the subset
of points of $X$ where all elements of all $R_j$ with
$j \in \{j_1,\ldots,j_p\}$ vanish has codimension not less than $p$
in $V$. Then $R$ has a system of parameters of degrees $d_1,\ldots,d_m$.
\qed
\end{proposition}

When applying this criterion it is convenient to have a notation
for `the codimension of the subset of $X$ defined by the vanishing
of all invariants with degree in $\{j_1,\ldots,j_p\}$'.
We'll use $[j_1,\ldots,j_p]$.

Note that for $e \ge 1$ an invariant $g^e$ vanishes if and only if
$g$ vanishes. It follows that if $j_h | j_i$, $h \ne i$, then
$[j_1,\ldots,j_p] = [j_1,\ldots,j_{h-1},j_{h+1},\ldots,j_p]$.

\subsection{From hsop degrees to generator degrees}
Let $R := \oo(V)^{\SLtwo}$. This is a graded algebra, and formulas
to compute its Poincar\'e series $P(t)$ are well-known.
If $(P_1,\ldots,P_m)$ is a system of parameters of $R$, with degree
sequence $(d_1,\ldots,d_m)$, then $P(t)$ can be
written as
$$P(t) = \frac{t^{e_1} + \ldots + t^{e_s}}{(1-t^{d_1}) \ldots (1-t^{d_m})}$$
and there exist homogeneous $G_1,\ldots,G_s \in R$ with degrees
$e_1,\ldots,e_s$, such that
\[
R = \notsobigoplus_{i=1}^s G_i \Cc[P_1,\ldots,P_m] .
\] 
Now $\{P_1,\ldots,P_m,G_1,\ldots,G_s\}$ is a (not necessarily minimal)
system of generators of $R$, and
$\max \{ d_1, \ldots , d_m, e_1, \ldots , e_s \}$ is an upper bound
for the degrees of a set of generators of $R$.

\subsection{Polarization}
Let $j = j(f)$ be an invariant of degree $d$ defined on forms $f \in V_n$.
Let ${\bf i} = (i_1,\ldots,i_s)$ be a sequence of nonnegative integers
with $\sum i_h = d$. The ${\bf i}$-{\em polarizations} $j_{\bf i}$ of $j$
are defined on $sV_n$ by
$$
j(\sum_i \lambda_i f_i) = \sum_{\bf i} j_{\bf i}(f_1,\ldots,f_s)
\lambda_1^{i_1} \ldots \lambda_s^{i_s} .
$$
Kraft \& Wallach \cite{KWpolar} showed for $n > 1$ that $\nn (sV_n)$
is defined by the polarizations of any set of functions defining $\nn (V_n)$.

\subsection{Transvectants}
The simplest examples of invariants are obtained using
{\em transvectants}. Given $g \in V_m$ and $h \in V_n$
the expression
$$
(g,h) \mapsto (g,h)_p:=\frac{(m-p)!(n-p)!}{m!n!}
\sum_{i=0}^p (-1)^i \binom{p}{i}
\frac{\partial ^p g}{\partial x^{p-i}\partial y^i}
\frac{\partial ^p h}{\partial x^i \partial y^{p-i}}
$$
defines a linear and $\SLtwo$-equivariant map
$V_m\otimes V_n\rightarrow V_{m+n-2p}$, which is classically called
the {\it p-th transvectant} (Ueberschiebung).
% scorrimento
The $(g,h)_p$ are the coefficients of the image of $g \otimes h$
under the isomorphism (Clebsch-Gordan formula)
\[V_m \otimes V_n \iso V_{m+n}\oplus V_{m+n-2}\oplus \ldots \oplus V_{m-n}\]
(for $m\geq n$).
We have $(f,g)_0=fg$ and $(f,f)_{2i+1}=0$ for all integers $i\geq 0$.

\section{Systems of parameters}\label{hsops}

In this section we give a homogeneous system of parameters
(or at least the degrees of a homogeneous system of parameters)
for the algebras of invariants in the cases occurring in Theorem \ref{main}.

% Finding generators in such cases was already routine in classical times
% (sometimes requiring months of work). The sets found classically were
% not necessarily minimal, but 
% 
% but to our knowledge this is the first time
% that explicit systems of parameters, or the degrees of such sets,
% were found in these cases.\footnote{But we do not give them.}
% So far even the minimal degrees of a system of parameters were unknown.
% 

\medskip\noindent
Let $\discr(f)$ denote the discriminant of the polynomial $f$.
If $f$ has degree $m$, then $\discr(f)$ is an expression of degree
$2m-2$ in the coefficients of $f$. This expression vanishes if and only if
$f$ has a root of multiplicity greater than 1.
% in terms of the roots: a_0^{2n-2} \prod_{i<j} (z_i-z_j)^2
% for f = ax+b this is 1

\medskip\noindent
Let $\res(f,g)$ denote the resultant of the polynomials $f$ and $g$.
If $f$ and $g$ have degrees $m$ and $n$, respectively, then
$\res(f,g)$ is an expression of degree $m+n$ in the coefficients
of $f$ and $g$. It vanishes if and only if $f$ and $g$ have a common root.

\medskip\noindent
Let $\sim$ denote equality up to a nonzero constant.

\begin{lemma}\label{one}\quad

(i) Let $l,m \in V_1$. Then $\res(l,m) = (l,m)_1$.

(ii) Let $q \in V_2$. Then $\discr(q) \sim (q,q)_2$.

(iii) Let $l \in V_1$ and $q \in V_2$. Then $\res(l,q) = (q,l^2)_2$.

(iv) Let $q,r \in V_2$. Then
$\res(q,r)  = (q,r)_2\supr[2] - (q,q)_2\,(r,r)_2$.

(v) Let $f \in V_3$. Then $\discr(f) \sim (f,(f,(f,f)_2)_1)_3$.

(vi) Let $l \in V_1$ and $f \in V_3$. Then $\res(l,f) \sim (f,l^3)_3$.
\qed
\end{lemma}

\noindent
Brion~\cite{Brion} shows for
$V = V_{\bf n} = V_{n_1} \oplus \ldots \oplus V_{n_p}$
with $p > 1$ that a multihomogeneous system of parameters exists in
only 13 cases, namely
% the cases
% ${\bf n} = (1,1)$, $(1,2)$, $(1,3)$, $(1,4)$, $(2,2)$, $(2,3)$, $(2,4)$,
% $(3,3)$, $(4,4)$, $(1,1,1)$, $(1,1,2)$, $(1,2,2)$, $(2,2,2)$,
precisely the cases with $p \in \{2,3\}$ in Popov's classification
(Table \ref{tabpopov}). We give the systems here (improving that for $2V_3$),
together with those for $V_2$, $V_3$, $V_4$, for later use.
% (For $2V_3$ Brion gives a hsop with degrees $4,4,4,4,4$. We give one
% with degrees $2,4,4,4,4$.)

\begin{proposition}\label{brion}
The modules $V$ given in the table below, with generic elements as indicated,
have a (multi-)homogeneous system of parameters as given.

\medskip\noindent\begin{tabular}{|@{~}c|c|@{~}c@{~}|l@{~}|}
\hline
$V$ & {element} & {hsop degrees} & ${\rm hsop}(V)$ \\
\hline
$V_2$ & $q$ & $2$ & $(q,q)_2$ \\
$V_3$ & $c$ & $4$ & $\discr(c)$ \\
$V_4$ & $f$ & $2,3$ & $(f,f)_4$, $(f,(f,f)_2)_4$ \\
\hline
$2V_1$ & $(l,m)$ & $2$ & $(l,m)_1$ \\
$V_1 \oplus V_2$ & $(l,q)$ & $2,3$ & $(q,q)_2$, $(q,l^2)_2$ \\
$V_1 \oplus V_3$ & $(l,c)$ & $4,4,4$ & ${\rm hsop}(V_3)$, $(c,l^3)_3$,
 $(c,(c,l^2)_1)_3$ \\
$V_1 \oplus V_4$ & $(l,f)$ & $2,3,5,6$ & ${\rm hsop}(V_4)$,
 $(f,l^4)_4$, $((f,f)_2,l^4)_4$ \\
$2V_2$ & $(q,r)$ & $2,2,2$ & polarized ${\rm hsop}(V_2)$ \\
$V_2 \oplus V_3$ & $(q,c)$ & $2,3,4,5$ & $(q,q)_2$,
 $(c,(c,q)_1)_3$, ${\rm hsop}(V_3)$, $\res(q,c)$ \\
$V_2 \oplus V_4$ & $(q,f)$ & $2,2,3,3,4$ &
 $(q,q)_2$, ${\rm hsop}(V_4)$, $(f,q^2)_4$, $((f,f)_2,q^2)_4$ \\
$2V_3$ & $(c,d)$ & $2,4,4,4,4$ & $(c,d)_3$, $\discr(c)$, $\discr(d)$, \\
 &&& $(c,(c,(c,d)_2)_1)_3$, $(d,(d,(c,d)_2)_1)_3$ \\
$2V_4$ & $(f,g)$ & $2,2,2,3,3,3,3$ &
 polarized ${\rm hsop}(V_4)$ \\
\hline
$3V_1$ & $(l,m,n)$ & $2,2,2$ & $(l,m)_1$, $(l,n)_1$, $(m,n)_1$ \\
$2V_1 \oplus V_2$ & $(l,m,q)$ & $2,2,3,3$ & $(l,m)_1$, $(q,q)_2$,
 $(q,l^2)_2$, $(q,m^2)_2$ \\
$V_1 \oplus 2V_2$ & $(l,q,r)$ & $2,2,2,3,3$ & ${\rm hsop}(2V_2)$,
 $(q,l^2)_2$, $(r,l^2)_2$\\
$3V_2$ & $(q,r,s)$ & $2,2,2,2,2,2$ & polarized ${\rm hsop}(V_2)$ \\
\hline
\end{tabular}
\end{proposition}

\subsection{The cases $V = pV_1$ and $V = pV_2$}
Let $V = pV_1$ (or $V = pV_2$), where $p > 1$.
Let $I$ be the ideal of $R$ generated by the invariants of degree 2.
% by the transvectants $(l_i,l_j)_2$ $(1 \le i < j \le p)$ resp.
% by the transvectants $(q_i,q_j)_2$ $(1 \le i \le j \le p)$.
By Lemma \ref{one} (i) (or (ii),(iv)) we have $\nn(V)=\vv(I)$.
Consider the proof of the Noether Normalization Lemma.
In the situation of an algebra where all generators are homogeneous
of the same degree, it produces a homogeneous set of parameters of
this same degree. Therefore, $I$, and hence $R$, has a homogeneous
system of parameters consisting of $2p-3$ (or $3p-3$) elements of degree 2.
(For generators, see \S\ref{onetwo} below.)

\subsection{The case $V = mV_1 \oplus nV_2$}
\begin{proposition}
Let $V = mV_1 \oplus nV_2$ with $m,n>0$, and let $R = \oo(V)^{\SLtwo}$.

If $2n+1 \ge m$, then $R$ has a system of parameters consisting
of $m+2n-2$ invariants of degree $2$ and $m+n-1$ invariants
of degree $3$.

If $2n+1 < m$, then $R$ has a system of parameters consisting
of $m+2n-2$ invariants of degree $2$ and $3n$ invariants
of degree $3$, and $m-2n-1$ invariants of degree $6$.
\end{proposition}
\begin{proof}
Invoke Dixmier's Criterion.
We have $\dim X = 2m+3n-3$. We have to show that $[2] \ge m+2n-2$
and $[3] \ge m+n-1$ (if $2n+1 \ge m$) and $[2,3] \ge 2m+3n-3$.

Consider $v = (l_1,\ldots,l_m,q_1,\ldots,q_n) \in V$.
According to Lemma \ref{one}, among the invariants of degree 2
there are $\res(l_i,l_j)$ and $\discr(q_i)$ and $\res(q_i,q_j)$
for all $i,j$ $(j \ne i)$, and among the invariants of degree 3
there are the $\res(l_i,q_j)$. If all of these vanish on $v$,
then $v \in \nn (V)$. Hence $[2,3] = \dim X$, as desired.

Look at $[2]$.
The element
$g = \left(\begin{smallmatrix} a & b \\ c & d \end{smallmatrix}\right)$
of $\SLtwo$ acts via $g.x = dx-by$, $g.y = -cx+ay$.
Using $\SLtwo$ one can move the first nonzero linear form (if there is one)
to $x$. Given that all linear forms have a common zero, this means that
the remaining at most $m-1$ linear forms now look like $cx$
for various constants $c$, zero or not. Given that all quadratic forms
have a common double root, and that elements
$\left(\begin{smallmatrix} 1 & 0 \\ c & 1 \end{smallmatrix}\right)$
preserve $x$ and act on the $q_i$, we can pick orbit representatives
$dy^2$ (for nonzero $d$) for the quadratic forms, unless none of them
involve $y$, in which case we are in the nullcone.
Altogether the result has dimension at most $m+n-1$, so that
$[2] \ge \dim X - (m+n-1) = m+2n-2$, as desired.

Look at $[3]$.
Now we are in a part of $X$ where each of the quadratic forms
shares a zero with each of the linear forms.
Suppose first that at least one linear form and at least one quadratic form
are nonzero. We can pick $x$ as orbit-representative for one linear form.
All quadratic forms now look like $x(ax+by)$ for various $a,b$,
and we can normalize one to $cx^2$ or $cxy$ with nonzero $c$.
All linear forms now look like $dx$ or $dy$ for various constants $d$.
Altogether the result has dimension at most $m-1+2(n-1)+1 = m+2n-2$,
as desired.

If all linear forms are zero, then consider the invariants
$(q_i,(q_j,q_k)_1)_2$ with $1 \le i < j < k \le n$ of degree 3.
Let us compute this for the case of the three forms
$ax^2+2bxy+cy^2$, $dx^2+2exy+fy^2$, $gx^2+2hxy+iy^2$.
The result is (up to a constant) $aei-afh-bdi+bfg+cdh-ceg$,
the determinant
$$\left|\begin{array}{ccc} a & b & c \\ d & e & f \\ g & h & i \end{array}
\right| .$$
If all such determinants vanish, then any three of the quadratic forms
are linearly dependent, so that $n$ quadratic forms involve at most
$2n+2-3 = 2n-1$ constants, after dividing out $\SLtwo$.
Since $m \ge 1$, this is not more than $m+2n-2$ and we also have
the desired bound in this case.

If all quadratic forms are zero, then the result has dimension $2m-3$,
so codimension $3n$. We showed $[3] \ge \min (m+n-1, 3n)$.
It follows that for $mV_1 \oplus nV_2$ with $2n \ge m-1 \ge 0$ there is a hsop
with $m+2n-2$ elements of degree 2 and $m+n-1$ elements of degree 3.
This proves the first claim.

For the second claim we still have to show that $[6] \ge 1$ and
$[2,6] \ge m+2n-1$ and $[3,6] \ge m+n-1$, but $[6] \ge [2,3] = 2m+3n-3$.
\end{proof}

\subsection{The case $V = 2V_1 \oplus V_3$}
Let $V = 2V_1 \oplus V_3$. We have $\dim X = 5$.
The ring $R$ has a homogeneous system of parameters
with degrees $2,4,4,4,4$.

Indeed, pick $(l,m,c) \in V$. By Proposition \ref{brion} we find
$(l,m) \in \nn (2V_1)$ and $(l,c),(m,c) \in \nn (V_1 \oplus V_3)$
(and hence $(l,m,c) \in \nn (V)$)
when all invariants of degree 4 vanish on $(l,m,c)$. This shows
that there is a hsop with degrees $4,4,4,4,4$. Since $[2] \ge 1$
there is also a hsop with degrees $2,4,4,4,4$ by Dixmier's Criterion.

\subsection{The case $V = V_1 \oplus V_2 \oplus V_3$}
Let $V = V_1 \oplus V_2 \oplus V_3$. We have $\dim X = 6$.
We show that the ring $R$ has a system of parameters with degrees
$3,3,4,4,4,5$.

A form on which all invariants of degrees $3,4,5$ vanish
is in the nullcone, so that $[3,4,5] = 6$. We have to check that
$[5] \ge 1$, $[3] \ge 2$, $[4] \ge 3$, $[3,5] \ge 3$,
$[4,5] \ge 4$, $[3,4] \ge 5$.
Let the forms be $l,q,f$.

If all invariants of degree 4 vanish, then $f$ has a double root
that is also a root of $l$, and $q$ has a double root and only
3 variables are left, so $[4] \ge 3$. If moreover $\res(q,f)$
vanishes, then only pieces of dimension 2 are left, so $[4,5] \ge 4$.
Or, if moreover $(q,l^2)_2$ and $(f,(f,q)_1)_3$ vanish,
then we are in the nullcone or in a piece of dimension at most 1,
so that $[3,4] \ge 5$.

If all invariants of degree 3 vanish (and in particular
$\res(q,l)$ and $(f,lq)_3$ and $(f,(f,q)_1)_3$), and $q \ne 0$,
then w.l.o.g. either
$q = x^2$, $l = ax$, $f = bx^3+3cx^2y+3dxy^2+ey^3$ with $ae=0$ and $ce-d^2=0$,
or $q = hxy$, $l = ax$, $f = bx^3+3cx^2y+3dxy^2+ey^3$ with $ae=0$ and $be=cd$,
of dimension at most 3 (since the torus of $\SLtwo$ still acts).
If $q = 0$ then we are in $V_1 \oplus V_3$ of dimension 3.
Altogether $[3,5] \ge [3] \ge 3$.
(In fact $[3,5] = 3$ because of the part with $q=0$.)

Finally $[5] \ge 1$ is clear.

Note that in this case the denominator of $P(t)$ might suggest
to look for a hsop with degrees $2,3,3,4,4,5$, but there is none
since $[2,3,5] = 3$.

\subsection{The case $V = V_1 \oplus V_2 \oplus V_4$}

Let $V = V_1 \oplus V_2 \oplus V_4$. We have $\dim X = 7$.
We show that the ring $R$ has a system of parameters with degrees
$2,2,3,3,4,5,6$.

Since $(l,q,f) \in \nn(V)$ if and only if $(l,q) \in \nn(V_1\oplus V_2)$
and $(l,f) \in \nn(V_1\oplus V_4)$
and $(q,f) \in \nn(V_2\oplus V_4)$, it follows that if
the eight invariants $(q,q)_2$, $(f,f)_4$, $(f,(f,f)_2)_4$,
$(q,l^2)_2$, $(f,q^2)_4$, $((f,f)_2,q^2)_4$, $(f,l^4)_4$,
$((f,f)_2,l^4)_4$ (of degrees $2,2,3,3,3,4,5,6$)
vanish, then $(l,q,f)\in \nn(V)$.

The above five invariants of degrees $2,4,5,6$, together
with two random linear combinations of the invariants of degree 3 will
constitute a hsop. In particular, we find that the two combinations
$(f,(f,f)_2)_4+(f,q^2)_4$ and $(q,l^2)_2 - (f,(f,f)_2)_4$ yield such a hsop.
(Using Singular one finds that the ideal generated by these seven invariants
contains the sixth power of each invariant of degree 3. Now apply
Proposition \ref{hilb}.)

Note that in this case the denominator of $P(t)$ suggests the
existence of a hsop with degrees $2,2,3,3,3,4,5$. But such a hsop does
not exist: all invariants of degrees 2 up to 5 vanish on
$(x,0,4xy^3)$, which is not in $\nn (V)$.

\subsection{The case $V = 2V_2 \oplus V_3$}
Let $V = 2V_2 \oplus V_3$. We have $\dim X = 7$.
The ring $R$ has a system of parameters with degrees
$2,2,3,4,5,5,6$.

Indeed, we have to check that $[3] \ge 1$, $[2] \ge 2$, $[5] \ge 2$,
$[4] \ge 3$, $[2,3] \ge 3$, $[3,5] \ge 3$, $[6] \ge 4$, $[2,5] \ge 4$,
$[3,4] \ge 4$, $[4,5] \ge 5$, $[4,6] \ge 5$, $[2,3,5] \ge 5$,
$[5,6] \ge 6$, $[3,4,5] \ge 6$, $[4,5,6] \ge 7$.

Let the forms be $q$, $r$, $f$.
If $f \ne 0$, then we can normalize $f$ to one of $x^3$, $x^2y$ or $axy(x+y)$.
If $f = 0$ but $q \ne 0$, then we can normalize $q$ to $x^2$ or $xy$.

If all invariants of degree 2 vanish, then $\discr(q)=\discr(r)=\res(q,r)=0$
so that $q$ and $r$ have a common double zero. Now if $q \ne 0$, then $r$
is determined by a single constant. So dimensions are at most 1 larger
than for the corresponding case for $V_2 \oplus V_3$.
Hence $[2] \ge 3$, $[4] \ge 4$, $[6] \ge [2,3] \ge 4$, $[2,5] \ge 4$,
$[4,6] \ge [3,4] \ge 5$, $[4,5] \ge 5$, $[2,3,5] \ge 5$, $[3,4,5] \ge 6$.

If all invariants of degree 3 vanish, then
$(f,(f,q)_1)_3 = (f,(f,r)_1)_3 = 0$.
That $(f,(f,q)_1)_3 = 0$ says that the three quadratic forms
$\frac{df}{dx}$, $\frac{df}{dy}$ and $q$ are linearly dependent.
So, here either $f$ has a triple root,
or $f_x$ and $f_y$ are linearly independent and $q,r$ are
determined by two coefficients each.
Hence $[3] \ge 2$.

If all invariants of degree 5 vanish, then $\res(f,q) = \res(f,r) = 0$
restrict $q,r$ when $f \ne 0$. Hence $[5] \ge 2$.

If all invariants of degrees $3,5$ vanish, then if $f$ does not have a
double root, then $q,r$ are determined by one coefficient each.
If $f = x^2y$ or $f = x^3$, then $q,r$ are determined by two coefficients each.
If $f=0$, then we are in $2V_2$ which has dimension 3.
Hence $[3,5] \ge 3$.

If all invariants of degrees $2,3,5$ vanish, then either $f = axy(x+y)$,
$q = r = 0$, or ($f = x^2y$ or $f = x^3$), $q = bx^2$, $r = cx^2$ and
$(q,r,f) \in \nn (V)$, or $f = 0$ and again $(q,r,f) \in \nn (V)$.
Hence $[5,6] \ge [2,3,5] \ge 6$ and $[3,4,5] \ge [2,3,5] \ge 6$.

Finally, if all invariants of degrees $2,3,4,5$ vanish, then $\discr(f) = 0$
and we are in the nullcone.

\medskip
Note that in this case the denominator of $P(t)$ might suggest
to look for a hsop with degrees $2,2,3,3,4,5,5$, but there is none
since $[3,5] = 3$.
% Look at the part with $f = x^2y$.
% Maybe conclude that at least 4 degrees must be even.

\subsection{The case $V = 2V_2 \oplus V_4$}

Let $V = 2V_2 \oplus V_4$. We have $\dim X = 8$.
We show that the ring $R$ has a system of parameters with degrees
$2,2,2,3,3,3,4,4$.

Since $(q,r,f) \in \nn(V)$ if and only if $(q,r) \in \nn(2V_2)$
and $(q,f),(r,f) \in \nn(V_2\oplus V_4)$,
it follows that if the nine invariants
$(q,q)_2$, $(q,r)_2$, $(r,r)_2$, $(f,f)_4$, $(f,(f,f)_2)_4$,
$(f,q^2)_4$, $(f,r^2)_4$, $((f,f)_2,q^2)_4$, and $((f,f)_2,r^2)_4$
(of degrees $2,2,2,\maysplit 2,3,3,3,4,4$) vanish,
then $(q,r,f)\in \nn (V)$.

The above five invariants of degrees $3,4$, together with
three random linear combinations of the four invariants of
degree 2 will constitute a hsop. In particular, we find that
the three combinations $(q,q)_2+(q,r)_2$, $(r,r)_2+(f,f)_4$,
and $(q,q)_2-(f,f)_4$ yield such a hsop.
(Using Singular one finds that the ideal generated by these eight invariants
contains the 7th power of each invariant of degree 2. Now apply
Proposition \ref{hilb}.)

Note that in this case the denominator of $P(t)$ might suggest
to look for a hsop with degrees $2,2,2,2,3,3,3,4$, but no such hsop
exists since all invariants of degree 2 or 3 vanish on
$(ax^2,bx^2,xy^3)$, so that $[2,3] \le 6$.

\subsection{The case $V = 3V_4$}

Let $V = 3V_4$. We have $\dim X = 12$.
We show that the ring $R$ has a system of parameters consisting of
6 invariants of degree 2 and 6 of degree 3.

Indeed, since $2V_4$ has a hsop consisting of elements
of degrees 2 and 3 only, $\nn (3V_4)$ is determined by
elements of degrees $2,3$, so that $[2,3] = 12$.
We have $\dim R_2 = 6$, and using Singular we find that $[2] \ge 6$.
Since $[2,3] = 12$, and the six invariants of degree 2 can decrease
dimensions by not more than 6, we must have $[3] \ge 6$.
By Dixmier's Criterion there is a hsop as claimed.

On the other hand, using Singular we also find that $[3] \ge 7$.
It follows that there is also a hsop consisting of 5 invariants
of degree 2 and 7 of degree 3.

\subsection{The case $V = 2V_1 + V_4$}
Let $V = 2V_1 + V_4$. We have $\dim X = 6$.
We show that the ring $R$ has a system of parameters with degrees
$2,3,5,5,6,6$.

Indeed, let the forms be $l,m,f$. The vanishing of
$(f,f)_4$, $(f,(f,f)_2)_4$, $(f,l^4)_4$, $((f,f)_2,l^4)_4$,
$(f,m^4)_4$, $((f,f)_2,m^4)_4$ implies that $l$ and $f$, and also $m$ and $f$,
have a common zero of multiplicity at least 3 for $f$.
If $f \ne 0$ then $(l,m,f) \in \nn (V)$. If $f = 0$ then the same conclusion
follows if also $(l,m)_1 = 0$.
By Dixmier's Criterion there is a hsop with the above degrees.
%
% We need $[2] \ge 1$, $[3] \ge 1$, $[5] \ge 2$, $[6] \ge 4$,
% $[2,3] \ge 2$, $[2,5] \ge 3$, $[3,5] \ge 3$, $[5,6] \ge 6$,
% $[2,3,5] \ge 4$. In the part with $f \ne 0$ or $l = 0$ or $m = 0$
% these inequalities hold. In the part with $f = 0$, $l \ne 0$, $m \ne 0$
% we have w.l.o.g. $l=x$, $m=ay$. Since that is only 1-dimensional
% the only thing to check is that $(l,m)_1 = 0$ when all invariants
% of degrees $2,3,4,5,6$ vanish.
% Probably one can give an explicit hsop by adding $(l,m)_1^3$ to
% $((f,f)_2,m^4)_4$.

\subsection{The case $V = V_3 + V_4$}

Let $V = V_3 + V_4$. We have $\dim X = 6$.
We show that the ring $R$ has a system of parameters with degrees
$2,3,4,5,6,7$.

Indeed, let $(f,g) \in V$. Let $i_2 = (g,g)_4$, $i_3 = (g,(g,g)_2)_4$,
$i_4 = \discr(f)$. The vanishing of $i_2,i_3,i_4$ forces $f$ to have
a double zero, and $g$ to have a triple zero.
If these zeros occur at different places,
then w.l.o.g. $f = x^2(ax+by)$ and $g = y^3(cx+dy)$.
Consider the four invariants
$i_5' = (f,(f,(f,(f,g)_1)_3)_1)_3$, $i_5'' = (f,(f,(g,(g,g)_2)_1)_3)_3$,
$i_6 = (f,(f,(f,(f,(g,g)_2)_1)_3)_1)_3$,
$i_7 = (f,(f,g)_3\supr[3])_3$.
They have values $b^4d$,  $a^2c^3$, $b^4c^2$ and
$ab^3c^3 - 10a^2b^2c^2d + 32a^3bcd^2 - 32a^4d^3$,
all up to some nonzero constant. Let $i_5 = i_5' + i_5''$.
% (an arbitrary linear combination with nonzero coefficients).
Suppose $i_5,i_6,i_7$ vanish on $(f,g)$.
If the form $f$ has a triple zero (i.e., $b=0$) then $i_5,i_6,i_7$
take the values $a^2c^3$, 0 and $a^4d^3$, so that either $f=0$ or $g=0$.
Otherwise we have w.l.o.g. $b=1$, the vanishing of $i_6$ and $i_5$ forces
$c=0$ and $d=0$, so that $g=0$.
This shows that $i_2,i_3,i_4,i_5,i_6,i_7$ determine the nullcone
and hence form a hsop.

\section{Generators}\label{gens}

In this section we give a set of generators for the algebras of invariants
in the cases occurring in Theorem \ref{main}.

\subsection{The case $V = V_1 \oplus V_2 \oplus V_3$}
Let $V = V_1 \oplus V_2 \oplus V_3$.
Here $r = 15$, cf. \cite{GordanS}, \cite{Sylv}, \cite{GrYo}, \S140.
Let the forms be $l$, $q$, $c$ of degrees 1, 2, 3, respectively.
The table below gives the 15 basic generators with degree and multidegree.
Abbreviations are as above.
$$
\medskip
\begin{tabular}{|@{~}c@{~}c@{~}l@{~}|@{~}c@{~}c@{~}l@{~}|@{~}c@{~}c@{~}l@{~}|}
\hline
dg & mdeg & ~~~form &
dg & mdeg & ~~~form &
dg & mdeg & ~~~form \\
\hline
2 & 020 & $(q,q)_2$ &
4 & 121 & $({\rm lq}_1,{\rm qc}_2)_1$ &
5 & 212 & $({\rm lq}_1,(l,{\rm cc}_2)_1)_1$ \\
3 & 012 & $(q,{\rm cc}_2)_2$ &
4 & 202 & $({\rm lc}_1,{\rm lc}_1)_2$ &
5 & 311 & $(l^2,({\rm lc}_1,q)_1)_2$ \\
3 & 111 & $(q,{\rm lc}_1)_2$ &
4 & 301 & $(l^2,{\rm lc}_1)_2$ &
6 & 123 & $((l,{\rm cc}_2)_1,(q,{\rm qc}_1)_2)_1$ \\
3 & 210 & $(l^2,q)_2$ &
5 & 032 & $({\rm qc}_2, (q,{\rm qc}_1)_2)_1$ &
6 & 303 & $((l,{\rm cc}_2)_1,(l^2,c)_2)_1$ \\
4 & 004 & $({\rm cc}_2,{\rm cc}_2)_2$ &
5 & 113 & $(({\rm lc}_1,q)_1,{\rm cc}_2)_2$ &
7 & 034 & $({\rm qc}_2^{\,2},({\rm qc}_1,c)_2)_2$ \\
\hline
\end{tabular}
$$
Given the Poincar\'e series and the hsop degrees, it suffices to check
that these generators generate $R_d$ for $d \le 14$, and they do.

\subsection{The case $V = V_1 \oplus V_2 \oplus V_4$}
Let $V = V_1 \oplus V_2 \oplus V_4$.
Here $r = 18$, cf. \cite{GordanS}, \cite{Sylv}.
Let the forms be $l$, $q$, $f$ of degrees 1, 2, 4, respectively.
The table below gives the 18 basic generators with degree and multidegree.
$$
\begin{tabular}{|@{~}c@{~}c@{~}l@{~}|@{~}c@{~}c@{~}l@{~}|@{~}c@{~}c@{~}l@{~}|}
\hline
dg & mdeg & ~~~form &
dg & mdeg & ~~~form &
dg & mdeg & ~~~form \\
\hline
2 & 020 & $(q,q)_2$ &
4 & 022 & $((f,f)_2,q^2)_4$ &
6 & 222 & $(f,l.(l,q)_1.(f,q)_2)_4$ \\
2 & 002 & $(f,f)_4$ &
5 & 401 & $(f,l^4)_4$ &
6 & 033 & $(f,(f,q)_2.(f,q^2)_3)_4$ \\
3 & 210 & $(q,l^2)_2$ &
5 & 221 & $(f,(lq.(l,q)_1)_4$ &
7 & 412 & $(f,lq.(f,l^3)_3)_4$ \\
3 & 021 & $(f,q^2)_4$ &
5 & 212 & $(f,l^2.(q,f)_2)_4$ &
7 & 223 & $(f,lq.(f,l.(f,q)_1)_4)_4$ \\
3 & 003 & $(f,(f,f)_2)_4$ &
6 & 411 & $(f,l^3.(l,q)_1)_4$ &
8 & 413 & $(f,l^2.(f,lq.(f,l)_1)_4)_4$\! \\
4 & 211 & $(f,l^2q)_4$ &
6 & 402 & $((f,f)_2,l^4)_4$ &
9 & 603 & $(f,l^3.(f,l.(f,l^2)_2)_3)_4$\! \\
\hline
\end{tabular}
$$
Given the Poincar\'e series and the hsop degrees, it suffices to check
that these generators generate $R_d$ for $d \le 15$, and they do.

\subsection{The case $V = 2V_2 \oplus V_3$}\label{V223}
Let $V = 2V_2 \oplus V_3$. We find $r = 18$.
Let the forms be $q$, $r$, $c$, of degrees 2, 2, 3, respectively.
Let $u = (c,q^2)_3$ and $v = (c,qr)_3$.
The table below gives the 18 basic generators with degree and multidegree.
(For multidegree $i.j.k$ only the entries with $i \ge j$ are given.)
$$
\begin{tabular}{|@{~}c@{~}c@{~}l@{~}|@{~}c@{~}c@{~}l@{~}|@{~}c@{~}c@{~}l@{~}|}
\hline
dg & mdeg & ~~~form &
dg & mdeg & ~~~form &
dg & mdeg & ~~~form \\
\hline
2 & 200 & $(q,q)_2$ &
4 & 004 & $(c,(c,(c,c)_2)_1)_3$ &
6 & 222 & $(c,(q,r)_1 v)_3$ \\
2 & 110 & $(q,r)_2$ &
5 & 302 & $(c,qu)_3$ &
7 & 304 & $(c,u.(c,(c,q)_1)_2)_3$ \\
3 & 102 & $(c,(c,q)_1)_3$ &
5 & 212 & $(c,qv)_3$ &
7 & 214 & $(c,u.(c,(c,r)_1)_2)_3$ \\
4 & 112 & $(c,q.(c,r)_2)_3$ &
6 & 312 & $(c,(q,r)_1 u)_3$ & && \\
\hline
\end{tabular}
$$
Given the Poincar\'e series and the hsop degrees, it suffices to check
that these generators generate $R_d$ for $d \le 17$, and they do.

\subsection{The case $V = 2V_2 \oplus V_4$}\label{V224}
Let $V = 2V_2 \oplus V_4$. We find $r = 19$.
Let the forms be $q$, $r$, $f$, of degrees 2, 2, 4, respectively.
The table below gives the 19 basic generators with degree and multidegree.
(For multidegree $i.j.k$ only the entries with $i \ge j$ are given.)
$$
\begin{tabular}{|c@{~~}c@{~}l|c@{~~}c@{~}l|c@{~~}c@{~~}l|}
\hline
dg & mdeg & ~~~form &
dg & mdeg & ~~~form &
dg & mdeg & ~~~form \\
\hline
2 & 200 & $(q,q)_2$ &
3 & 111 & $(f,qr)_4$ &
4 & 211 & $(f,q.(q,r)_1)_4$ \\
2 & 110 & $(q,r)_2$ &
3 & 003 & $(f,(f,f)_2)_4$ &
5 & 212 & $(f,(f,q)_2.(q,r)_1)_4$ \\
2 & 002 & $(f,f)_4$ &
4 & 202 & $(f,q.(f,q)_2)_4$ &
6 & 303 & $(f,(f,q)_2.(f,q^2)_3)_4$ \\
3 & 201 & $(f,q^2)_4$ &
4 & 112 & $(f,q.(f,r)_2)_4$ &
6 & 213 & $(f,(f,q)_2.(f,qr)_3)_4$ \\
\hline
\end{tabular}
$$
Given the Poincar\'e series and the hsop degrees, it suffices to check
that these generators generate $R_d$ for $d \le 12$, and they do.

\subsection{The case $V = 3V_4$}
Let $V = 3V_4$. We find $r = 25$.
Let the forms be $f$, $g$, $h$, all of degree 4. The table below
gives the 25 basic generators with degree and multidegree.
(For multidegree $i.j.k$ only the entries with $i \ge j \ge k$ are given.)
$$
\begin{tabular}{|@{~}c@{~}c@{~}l@{~}|@{~}c@{~}c@{~}l@{~}|@{~}c@{~}c@{~}l@{~}|}
\hline
dg & mdeg & ~~~form &
dg & mdeg & ~~~form &
dg & mdeg & ~~~form \\
\hline
2 & 200 & $(f,f)_4$ & 
3 & 210 & $(f,(f,g)_2)_4$ &
4 & 211 & $(f,(f,(g,h)_2)_2)_4$ \\
2 & 110 & $(f,g)_4$ & 
3 & 111 & $(f,(g,h)_2)_4$ &
5 & 221 & $(h,(f,g)_3.(f,g)_3)_4$ \\
3 & 300 & $(f,(f,f)_2)_4$ &
4 & 220 & $(f,(f,(g,g)_2)_2)_4$ & && \\
% 5 & 221 & $(f,(f,(g,(g,h)_1)_3)_2)_4$ \\
\hline
\end{tabular}
$$
Given the Poincar\'e series and the hsop degrees,
it suffices to check that these generators generate $R_d$ for
$d \le 15$, and they do. For $d = 15$ this required computing
the rank (34734) of a matrix with $10^{10}$ entries.
% matrix width needed: 288495
% 64:47:01elapsed
% For $pV_4$, cf.~\cite{Young}.

\subsection{The case $V = V_3 \oplus V_4$}
Let $V = V_3 \oplus V_4$. We find $r = 20$.
Let the forms be $c$, $e$.
Omit parentheses where that does not introduce ambiguity,
so that $(c,(c,e)_3(c,(c,e)_1)_3)_3$ is written $(c,ce_3 cce_{13})_3$.
Write $l = ce_3$.
The table below gives the 20 basic generators with degree and multidegree.
$$
\begin{tabular}{|@{~}c@{~}c@{~}l@{~}|@{~}c@{~}c@{~}l@{~}|@{~}c@{~}c@{~}l@{~}|}
\hline
dg & mdeg & ~~~form &
dg & mdeg & ~~~form &
dg & mdeg & ~~~form \\
\hline
2 & 02 & $ee_4$ &
7 & 43 & $cccceee_{213213}$ &
9 & 63 & $(c,cccce_{1231} cee_{23})_3$ \\
3 & 03 & $eee_{24}$ &
7 & 43 & $(c,cee_{23} cce_{13})_3$ &
9 & 63 & $(c,cccee_{2123} cce_{13})_3$ \\
4 & 40 & $cccc_{213}$ &
7 & 43 & $(c,l^3)_3$ &
9 & 63 & $(c,ccce_{123} l^2)_3$ \\
5 & 41 & $cccce_{1313}$ &
8 & 62 & $(c,cccce_{1231} l)_3$ &
9 & 45 & $(c,cee_{23} cee_{23} l)_3$ \\
5 & 23 & $cceee_{2133}$ &
8 & 44 & $(c,cceee_{2132} l)_3$ &
10 & 64 & $(c,cccee_{2123} ccee_{222})_3$ \\
6 & 42 & $ccccee_{21313}$ &
8 & 44 & $(c,cee_{23} l^2)_3$ &
10 & 64 & $(c,cccee_{2123} l^2)_3$ \\
6 & 42 & $(c,cce_{13} l)_3$ &
&&&
11 & 65 & $(c,ccceee_{21313} l^2)_3$ \\
\hline
\end{tabular}
$$
Given the Poincar\'e series and the hsop degrees, it suffices to check
that these generators generate $R_d$ for $d \le 18$, and they do.

\subsection{The case $V = W \oplus mV_1$}\label{mult1}
Given a basic system of covariants for a module $W$, one finds
a basic system of covariants for $W \oplus V_1$ by replacing
each covariant $j$ (of order $o$) of the system for $W$
by the $o+1$ covariants $(j,l^i)_i$ ($0 \le i \le o$)
and adding the covariant $l$, where $l$ is the linear form
corresponding to $V_1$. (This is classical. See \cite[\S 55]{Clebsch},
\cite[\S 138A]{GrYo}.)
Equivalently, given the invariants of $W \oplus V_1$, one finds
the invariants of $W \oplus mV_1$ by polarization.

One finds further results that should be regarded classical:

\begin{table}[h]
\centering
\begin{tabular}{lc|lc|lc}
module & $r$ & module & $r$ & module & $r$ \\
\hline
$2V_1 \oplus V_2$ & 5 &
$2V_1 \oplus V_3$ & 13 &
$2V_1 \oplus V_4$ & 20 \\
$2V_1 \oplus 2V_2$ & 13 &
$3V_1 \oplus 2V_2$ & 24 & % \cite{Perrin}
$3V_1 \oplus V_3$ & 30 \\
$2V_1 \oplus V_2 \oplus V_3$ & 35 &
$2V_1 \oplus 2V_4$ & 103 &
$3V_1 \oplus V_4$ & 63 \\
$2V_1 \oplus V_2 \oplus V_4$ & 57 & &&& \\
\end{tabular}
\caption{Some classical results involving multiple $V_1$}\label{clmult1}
\end{table}

\subsection{The case $V = mV_1 \oplus nV_2$}\label{onetwo}
%
% Algorithms for computing generators of the algebra of invariants of
% $pV_1$, $pV_2$, and $pV_1\oplus qV_2$ are described in \cite{GrYo}.
% Also, Gordan \cite{Gordan}, p.\,155 (and Kraft \& Weyman \cite{KrWe})
% shows that the algebra of invariants of $pV_2$ is generated
% by invariants of degrees 2 and 3.
%
Let $V = mV_1 + nV_2$. Then $r = \binom{n}{3} +
\binom{m+1}{2}\binom{n+1}{2} + \binom{m}{2} + \binom{n+1}{2}$.

If the forms are $\ell_i$ ($1 \le i \le m$) and $q_j$ ($1 \le j \le n$),
then there are $\binom{m}{2} + \binom{n+1}{2}$
basic invariants of degree 2, namely $(\ell _i,\ell _j)_1$ for $i<j$ and
$(q_i,q_j)_2$ for $i \le j$, and $n\binom{m+1}{2} + \binom{n}{3}$
basic invariants of degree 3, namely $(q_k,\ell _i \ell _j)_2$ for $i \le j$
and $(q_i,(q_j,q_k)_1)_2$ for $i < j < k$, and $\binom{m+1}{2}\binom{n}{2}$
basic invariants of degree 4, namely $((q_i,q_j)_1,\ell _k \ell _m)_2$
for $i < j$, $k \le m$.

In order to show this, we quote the following result
\cite[\S 54]{Clebsch}:

\begin{proposition}
Let $\mathcal R$ and $\mathcal S$ be two $\SLtwo $-algebras whose
covariants are finitely generated. Then the covariants of $\mathcal R
\oplus \mathcal S$ are alsofinitely generated. If $P_1,\ldots ,P_r$
are the generators of the covariants of $\mathcal R$, and $Q_1,\ldots
,Q_s$ are the generators of the covariants of $\mathcal S$, then a
finite generating system can be chosen from the set of transvectants
$[P,Q]_l$, $l\geq 0$, where $P$ is a monomial in the $P_i$'s and $Q$ a
monomial in the $Q_j$'s.
\end{proposition}

Apply this with ${\mathcal R} = mV_1$ and ${\mathcal S} = nV_2$,
with forms as above.
The covariants of $mV_1$ are generated by the $\ell_i$ themselves,
and the invariants $(\ell _i,\ell _j)_1$ for $i<j$.
The covariants of $nV_2$ are generated by the $q_i$ themselves,
the covariants $(q_i,q_j)_1$ for $i \le j$,
and the invariants $(q_i,q_j)_2$ for $i \le j$ and
$(q_i,(q_j,q_k)_1)_2$ for $i < j < k$.

We add to the set of generators of $R$ the invariants of degrees 3
$(q_k,\ell _i\ell _j)_2$ for $i \le j$, and the invariants of degree 4
$((q_i,q_j)_1,\ell _k \ell _m)_2$ for $i < j$, $k \le m$. Given that
$$ (r_1 \ldots r_p,\ell _1 \ldots \ell _{2p})_{2p}\sim \sum (r_1,\ell
_{i_1} \ell _{i_2})_2\ldots (r_p,\ell _{i_{2p-1}} \ell _{i_{2p}})_2,
$$ there are no other irreducible invariants.

\medskip\noindent
Addresses of authors:

\smallskip
\begin{minipage}{2in}
Andries E. Brouwer \\
Dept. of Math. \\
Techn. Univ. Eindhoven \\
P. O. Box 513 \\
5600MB Eindhoven \\
Netherlands \\
{\tt aeb@cwi.nl}
\end{minipage}
\begin{minipage}{2in}
Mihaela Popoviciu \\
Mathematisches Institut \\
Universit\"at Basel \\
Rheinsprung 21 \\
CH-4051 Basel \\
Switzerland \\
{\tt mihaela.popoviciu@unibas.ch}
\end{minipage}

\end{document}